\documentclass[natbib]{svcyclop}
\bibpunct{[}{]}{;}{n}{}{,} 
\usepackage{latexsym}
\usepackage{graphicx}
\usepackage{amsmath}
\usepackage{amssymb}

\newcommand{\ens}[1]{\mathbb{#1}}

\newcommand{\N}{\mathbb{N}}

\newcommand{\R}{\mathbb{R}}

\def\cal{\mathcal}

\def\N{{\cal N}}

\def\derpar#1#2{\frac{\partial#1}{\partial#2}}

\begin{document}

\title{Kinetic equations: computation}

\author{Lorenzo Pareschi\inst{1}} 

\institute{Dipartimento di Matematica, Universit\'a di Ferrara\\    
Via Machiavelli 35, 44121 Ferrara, Italy\\
E-mail: lorenzo.pareschi@unife.it}

\maketitle
\section{Synonyms}
Collisional equations, Transport equations, Boltzmann equations

\section{Mathematics Subject Classification}
65D32, 65M70, 65L04, 68Q25, 82C40

\section{Short Definition}
Kinetic equations bridge the gap between a microscopic description
and a macroscopic description of the physical reality. Due to the
high dimensionality the construction of numerical methods represents
a challenge and requires a careful balance between accuracy and 
computational complexity.

\section{Description}
\subsection{Kinetic equations}

Particle systems can be described at the microscopic level by
particle dynamics, i.e. systems of differential equations
describing the individual motions of the particles. However, 
they are extremely costly from a numerical
point of view, and bring little intuition on how a large
particle system behaves. Therefore, one is led to seek reduced
descriptions of particle systems which still preserve an
accurate description of the physical phenomena.
Kinetic models intend to describe particle systems by means of a
distribution function $f(x,v,t)$. This object represents a
number density in phase space, i.e. $f \, dx \, dv$ is the
number of particles in a small volume $dx \, dv$ in
position-velocity space about the point $(x,v)$ of this space.

In this short entry we will focus on computational methods for the interacting particle case described 
by the Boltzmann equation. This is motivated by its relevance for
applications and by the fact that it contains all major difficulties
present in other kinetic equations. From a numerical perspective, most of 
the difficulties are due to the multidimensional structure
of the distribution function. 
In particular the approximation of the 
collisional integral is a real challenge for numerical methods, since the 
integration runs on a highly-dimensional unflat manifold and is at
the basis of the macroscopic properties of the equation. Further
difficulties are represented by the presence of fluid-kinetic interfaces 
and multiple scales where most numerical methods loose their efficiency
 because they are forced to operate on a very short time scale.

Although here we review briefly only deterministic numerical methods let
us mention that several realistic numerical simulations are based on
Monte-Carlo techniques~\citep{Bi94, Na80, RjWa-96}. In the next paragraphs we summarize the main ideas at the basis
of two of the most popular way to approximate the distribution function in the velocity space, namely the discrete-velocity method~\citep{RS94,BoVaPaSc:cons:95,Bu:96,HePa:DVM:02} and the spectral method\cite{PP96,PaRu:spec:00,MP:03,FM:11,PTV}. Finally we shortly introduce the basic principles for the construction of schemes which are robust in fluid regions\cite{FJ:10,GaPaTo:97,DP:10}.

\subsection{Boltzmann equation}

Taking into account only binary interactions, the behavior of a
dilute gas of particles is described by the Boltzmann
equation~\cite{CIP95, Vill:hand}
 \begin{equation}
 \derpar{f}{t} + v \cdot \nabla_x f = Q(f,f)
 \end{equation}
where $f(t,x,v)$, $x,v \in \R^d$ ($d \ge 2$), is the
time-dependent particle distribution function in the phase space and
the collision operator $Q$ is defined by
 \begin{equation}\label{eq:Q}
 Q (f,f)(v) = \int_{v_* \in \R^d}
 \int_{\sigma \in \ens{S}^{d-1}}  B(\cos \theta,|v-v_*|) \,
 \left[ f'_* f' - f_* f \right] \, d\sigma \, dv_*.
 \end{equation}
Time and position act only as parameters in $Q$ and therefore
will be omitted in its description.
In~\eqref{eq:Q} we used the shorthands $f = f(v)$, $f_* = f(v_*)$,
$f ^{'} = f(v')$, $f_* ^{'} = f(v_* ^{'})$. The velocities of the
colliding pairs $(v,v_*)$ and $(v',v'_*)$ are related by
 \[
 v' = \frac{v+v_*}{2} + \frac{|v-v_*|}{2} \sigma, \qquad
 v'_* = \frac{v+v^*}{2} - \frac{|v-v_*|}{2} \sigma\nonumber.
 \]
The collision kernel $B$ is a non-negative function which only
depends on $|v-v_*|$ and $\cos \theta = ((v-v_*)/|v-v_*|)\cdot
\sigma$.
Boltzmann's collision operator has the fundamental properties of
conserving mass, momentum and energy
\begin{equation}
 \int_{v\in{\R}^d}Q(f,f) \, \phi(v)\,dv = 0. \qquad
 \phi(v)=1,v,|v|^2 
\end{equation}
Moreover, any equilibrium distribution
function $M$ such that $Q(M,M)=0$ has the form of a locally Maxwellian distribution
 \begin{equation}
 M(\rho,u,T)(v)=\frac{\rho}{(2\pi T)^{d/2}}
 \exp \left( - \frac{\vert u - v \vert^2} {2T} \right), 
 \end{equation}
where $\rho,\,u,\,T$ are the density, mean velocity
and temperature of the gas
 \begin{equation}
 \rho = \int_{v\in{\R}^d}f(v)dv, \quad u =
 \frac{1}{\rho}\int_{v\in{\R}^d}vf(v)dv, \quad T = {1\over{d\rho}}
 \int_{v\in{\R}^d}\vert u - v \vert^2f(v)dv.
 \end{equation}

\subsection{Discrete velocity methods}
Historically this was the first method for discretizing the
Boltzmann equation in velocity space. The discretization is built starting
from physical rather then numerical considerations. We assume the
gas particles can attain only a finite set of velocities \[
V_{\N}=\{v_1,v_2,v_3,\ldots,v_{\N}\},\quad v_i\in\R^d,
\]
and denote by $f_j(x,t)=f(v_j,x,t)$, $j=1,\ldots,\N$.
The collision pair
$(v_i,v_j)\leftrightarrow(v_k,v_l)$ is {admissible}
if {$v_i,v_j,v_k,v_l\,\in V_{\N}$} and preserves momentum and
energy {\[ v_i+v_j=v_k+v_l,\qquad
|v_i|^2+|v_j|^2=|v_k|^2+|v_l|^2.
\]}
The set of admissible output pairs {$(v_k,v_l)$}
corresponding to a given input pair {$(v_i,v_j)$} will be
denoted by {$C_{ij}$}.

{ The discrete collision operator is obtained as a quadrature formula based on the weights {$a_{ij}^{kl}$} related to the
collision {$(v_i,v_j)\leftrightarrow(v_k,v_l)$} which must
satisfy the relations {\[a_{ij}^{kl}\geq 0,\qquad
\sum_{k,l=1}^{\N} a_{ij}^{kl}=1,\,\forall i,j=1,\ldots,{\N}.
\]}
Next we introduce the {transition rates}
{$A_{ij}^{kl}=S|v_i-v_j|a_{ij}^{kl}$}, where {$S$} is
the cross sectional area of particles, and write the discrete
Boltzmann equation as {\[ \frac{\partial f_i}{\partial
t}+v_i\cdot \nabla_x f_i=Q_i(f,f),
\]}
with {\[ Q_i(f,f)=\sum_{{j,k,l=1}\atop{k,l\in C_{ij}}}^{\N} A_{ij}^{kl}(f_k f_l - f_i
f_j).
\]}
}
The discretized Boltzmann equation has the nice property of
preserving the essential physical features (conservations,
H-theorem, equilibrium states). However, from a computational point of view
the discrete Boltzmann equation presents two main
drawbacks. First the computational cost is larger then $O(\N^2)$ and second the accuracy is rather poor, less then first order accurate (see \cite{HePa:DVM:02} for example). 

\subsection{Spectral methods}
Spectral methods have been constructed recently with the goal to compensate the drawbacks of discrete velocity approximation. For the sake of simplicity we summarize their derivation in the case of the space homogeneous Boltzmann equations, although the schemes can be effectively used to compute the collision integral in a general setting.

The approximate function $f_N$ is represented as the truncated
Fourier series
 \[ f_N (v) = \sum_{k=-N}^{N} \hat{f}_k e^{i k \cdot v}, \
 \ \ \ \hat{f}_k = \frac{1}{(2 \pi)^d} \, \int _{\mathcal{D}_\pi} f(v) e^{-ik \cdot v} \, dv.\]
The spectral equation is the projection of the collision integral $Q^R(f,f)$, truncated over the ball of radius $R$ centered in the origin, in $\ens{P}^N$, the $(2N+1)^d$-dimensional
vector space of trigonometric polynomials of degree at most $N$ i.e.
 \begin{equation*}
 \derpar{f_N}{t} = \mathcal{P}_N Q^R (f_N,f_N)
 \end{equation*}
where $\mathcal{P}_N$ denotes the orthogonal projection on
$\ens{P}^N$ in $L^2(\mathcal{D}_\pi)$.
A straightforward
computation leads to the following set of ordinary differential
equations 
 \begin{equation}\label{eq:ode}
 \frac{d\hat{f}_k(t)}{dt} =
 \sum_{\underset{l+m=k}{l,m=-N}}^{N} \hat{\beta}(l,m) \, \hat{f}_l \, \hat{f}_m, \ \ \
 k=-N,...,N
 \end{equation}
where $\hat{\beta}(l,m)$ are the {\em kernel modes}, given by
$\hat{\beta}(l,m) = \beta(l,m) - \beta(m,m)$
with
 \begin{equation*}
 \beta(l,m) = \int_{x \in \mathcal{B}_R} \int_{y \in \mathcal{B}_R}
 \tilde{B}(x,y) \, \delta(x \cdot y) \,
 e^{i l \cdot x} \, e^{i m \cdot y} \, dx \, dy,
 \end{equation*}
 and
 \[ \tilde{B}(x,y) = 2^{d-1} \, B\left(-\frac{x \cdot (x+y)}{|x| |x+y|},|x+y|\right) \, |x+y|^{-(d-2)}. \]
As shown in\cite{MP:03} when $B$ satisfies the {\em decoupling assumption}
 $\tilde{B}(x,y) = a(|x|) \, b(|y|)$,
it is possible to approximate each $\hat{\beta}(l,m)$ by a sum
 \begin{equation}\label{eq:dec}
 \beta (l,m) \simeq \sum_{p=1} ^{A} \alpha_p (l) \alpha' _p (m).
 \end{equation}
This gives a sum of $A$ discrete convolutions, with $A\ll N$, and by standard FFT techniques a computational cost of $O(A\,N^d\log_2 N)$. Denoting by $\N=(2N+1)^d$ the total number of grid points this is equivalent to $O(A \, \N \log_2 \N)$ instead of $O(\N^2)$. Moreover, one gets the following consistency
result of spectral accuracy~\cite{MP:03}
 \begin{theorem}\label{theo:consist}
 For all $k > d-1$ such that $f \in H^k _p$
  \begin{equation*}
  \| Q^R (f,f) -  \mathcal{P}_N Q^{R,M} (f_N,f_N) \|_{L^2} \le C_1 \, \frac{ R^k \|f_N\|^2 _{H^k _p}}{M^k} +
                    \frac{C_2}{N^k} \, \Big( \|f\| _{H^k _p} + \|Q^R (f_N,f_N)\|_{H^k _p} \Big).
  \end{equation*}
 \end{theorem}

\subsection{Asymptotic-preserving methods}
Let us now consider the time discretization of the scaled Boltzmann equation 
\begin{equation}
 \derpar{f}{t} + v \cdot \nabla_x f = \frac{1}{\varepsilon}Q(f,f)
 \label{eq:beq}
 \end{equation}
where $\varepsilon>0$ is the Knudsen number. For small value of $\varepsilon$ we have a stiff problem and standard time discretization methods are forced to operate on a very small time scale. On the other hand in such regime formally $Q(f,f)\approx 0$ and the distribution function is close to a local Maxwellian. Thus the moments of the Boltzmann equation are well-approximated by the solution to the Euler equations of fluid-dynamics
\begin{equation}
\partial_t u+\nabla_x\cdot F(u)=0,
\label{eq:Euler} \end{equation} with\vskip -1.5cm
\[
u=(\rho,w,E)^T,\qquad F(u)=(\rho w, \varrho w \otimes (w+pI),
Ew+pw)^T,\quad p=\rho T,
\]
where $I$ is the identity matrix and $\otimes$ denotes the tensor
product.

We say that a time discretization method for (\ref{eq:beq}) of
stepsize $\Delta t$ is {\em asymptotic preserving (AP)} if,
independently of the
stepsize $\Delta t$, in the limit $\varepsilon\to 0$ becomes a
consistent time discretization method for the reduced system
(\ref{eq:Euler}).


When $\varepsilon \ll 1$ the problem is {\em stiff} and we must resort on implicit integrator to avoid small time step restriction. This however requires the inversion of the collision integral $Q(f,f)$ which is prohibitively expensive from the computational viewpoint.

On the other hand when $f \approx M[f]$ we know that the
collision operator $Q(f,f)$ is well approximated by its linear counterpart $Q(M,f)$ or by a simple relaxation
operator $(M-f)$. If we denote by $L(f)$ the selected approximate linear operator we can rewrite the equation introducing a penalization term as 
\begin{equation}
 \derpar{f}{t} + v \cdot \nabla_x f = \frac{1}{\varepsilon}{(Q(f,f)-L(f))}+\frac{1}{\varepsilon} L(f).
 \label{eq:beq2}
 \end{equation}
The idea now is to be implicit (or
exact) in the linear part $L(f)$ and explicit in the deviations
from equilibrium $Q(f,f)-L(f)$. This approach has been successfully presented in \cite{FJ:10} using implicit-explicit integrators and in \cite{GaPaTo:97,DP:10} by means of exponential techniques.

\section{Conclusions}
Computational methods for kinetic equations represent an emerging field in scientific computing. This is testified by the large amount of scientific papers which has been produced on the subject in recent years. We do not seek to review all of them here and focused our attention to the challenging case of the Boltzmann equation of rarefied gas dynamic. The major difficulties in this case are represented by the discretization of the multi-dimensional integral describing the collision process and by the presence of multiple time scales. Fast algorithms and robust stiff solvers are then essential ingredients of computational methods for kinetic equations.


\begin{thebibliography}{99}

\bibitem[{Bird(1994)}]{Bi94}
Bird G (1994) Molecular gas dynamics and direct simulation of gas flows.
  Clarendon Press, Oxford

\bibitem[{Bobylev et~al(1995)Bobylev, Palczewski, and
  Schneider}]{BoVaPaSc:cons:95}
Bobylev A, Palczewski A, Schneider J (1995) On approximation of the {B}oltzmann
  equation by discrete velocity models. C R Acad Sci Parais S\'er I Math
  320:639--644

\bibitem[{Buet(1996)}]{Bu:96}
Buet C (1996) A discrete velocity scheme for the {B}oltzmann operator of
  rarefied gas dynamics. Trans Theo Stat Phys 25:33--60

\bibitem[{Cercignani et~al(1994)Cercignani, Illner, and Pulvirenti}]{CIP95}
Cercignani C, Illner R, Pulvirenti M (1994) The mathematical theory of dilute
  gases, vol 106. Applied Mathematical Sciences

\bibitem[{Dimarco and Pareschi(2011)}]{DP:10}
Dimarco G, Pareschi L (2011) Exponential {R}unge-{K}utta methods for stiff
  kinetic equations. SIAM J Num Anal 49:2057--2077

\bibitem[{Filbet and Jin(2010)}]{FJ:10}
Filbet F, Jin S (2010) A class of asymptotic-preserving schemes for kinetic
  equations and related problems with stiff sources. J Comput Phys
  229:7625--7648

\bibitem[{Filbet and Mouhot(2011)}]{FM:11}
Filbet F, Mouhot C (2011) Analysis of spectral methods for the homogeneous
  {B}oltzmann equation. Trans Amer Math Soc 363:1947--1980

\bibitem[{Gabetta et~al(1997)Gabetta, Pareschi, and Toscani}]{GaPaTo:97}
Gabetta E, Pareschi L, Toscani G (1997) Relaxation schemes for nonlinear
  kinetic equations. SIAM J Numer Anal 34:2168--2194

\bibitem[{Mouhot and Pareschi(2006)}]{MP:03}
Mouhot C, Pareschi L (2006) Fast algorithms for computing the {B}oltzmann
  collision operator. Math Comp 75(256):1833--1852 (electronic)

\bibitem[{Nanbu(1980)}]{Na80}
Nanbu K (1980) Direct simulation scheme derived from the {B}oltzmann equation
  i. monocomponent gases. J Phys Soc Japan 49:2042--2049

\bibitem[{Panferov and Heintz(2002)}]{HePa:DVM:02}
Panferov V, Heintz A (2002) A new consistent discrete-velocity model for the
  {B}oltzmann equation. Math Methods Appl Sci 25:571--593

\bibitem[{Pareschi and Perthame(1996)}]{PP96}
Pareschi L, Perthame B (1996) A spectral method for the homogeneous {B}oltzmann
  equation. Trans Theo Stat Phys 25:369--383

\bibitem[{Pareschi and Russo(2000)}]{PaRu:spec:00}
Pareschi L, Russo G (2000) Numerical solution of the {B}oltzmann equation i.
  spectrally accurate approximation of the collision operator. SIAM J Numer
  Anal 37:1217--1245

\bibitem[{Pareschi et~al(2003)Pareschi, Toscani, and Villani}]{PTV}
Pareschi L, Toscani G, Villani C (2003) Spectral methods for the non cut-off
  {B}oltzmann equation and numerical grazing collision limit. Numerische
  Mathematik 93:527--548

\bibitem[{Rjasanow and Wagner(1996)}]{RjWa-96}
Rjasanow S, Wagner W (1996) A stochastic weighted particle method for the
  {B}oltzmann equation. J Comput Phys 124:243--253

\bibitem[{Rogier and Schneider(1994)}]{RS94}
Rogier F, Schneider J (1994) A direct method for solving the {B}oltzmann
  equation. Trans Theo Stat Phys 23:313--338

\bibitem[{Villani(2002)}]{Vill:hand}
Villani C (2002) A survey of mathematical topics in kinetic theory. Handbook of
  fluid mechanics, S. Friedlander and D. Serre, Eds. Elsevier Publ.

\end{thebibliography}
\end{document}